\documentclass[amsfont,amscd,11pt]{amsart}
\usepackage{amssymb}
\usepackage{amsmath}
\usepackage{amscd}
\usepackage{amssymb}
\usepackage{latexsym}
\newtheorem{theorem}{Theorem}[section]\newtheorem{lemma}[theorem]{Lemma}

\theoremstyle{definition}
\theoremstyle{remark}
\newtheorem{remark}[theorem]{Remark}\numberwithin{equation}{section}
\author{Bazanfar\'e Mahaman }
 \address{\hphantom{iii}Bazanfar\'e, D\'epartement de Math\'ematiques et Informatique \newline\hphantom{iii} Universit\'e Abdou Moumouni\newline\hphantom{iii} B.P. 10662, Niamey,Niger}
\email{bmahaman@yahoo.fr}

\title[Topology of manifolds with\ldots]{Topology of manifolds with asymptotically nonnegative Ricci curvature }

  \date{}
\thanks{2000 mathematics subject classification Primary 53C21, Secondary 53C20}
\begin{document}
\maketitle

\begin{abstract}
 In this paper, we study the topology of complete noncompact Riemannian manifolds with asymptotically nonnegative Ricci curvature. We show that a complete
noncompact manifold with asymptoticaly nonnegative Ricci curvature and sectional curvature $K_{M}(x)\ge -\frac{C}{d_{p}(x)^{\alpha}}$  is diffeomorphic to a Euclidean n-space $\mathbb{R}^n$ under some conditions on the density of rays starting from the base point $p$ or on the volume growth of geodesic balls in $M.$
\end{abstract}

\section{Introduction}
One of most important problems in Riemannian geometry is to find conditions under which manifold is of finite topological type: A manifold is said to have finite topological type if there exists a compact domain $\Omega$ with boundary such that $M\setminus\Omega$ is homeomorphic to $\partial\Omega\times[0,\infty[.$ 
 The fundamental notion involved in such a finite topological type
result is that of the critical point of a distance function introduced by Grove and
Shiohama \cite{8}.
Let $p$ a fix point and set $d_{p}(x) = d(p,x$)
A point $x\ne p$ is called critical point of $d_{p}$ if for any $v$ in the tangent space $T_{x}M$ there is minimal geodesic $\gamma$ from $x$ to $p$ forming an angle less or equal to $\pi/2$ with $\gamma'(0)$ (see \cite{8}).
 
In several papers it has been proved results for manifolds with nonnegative curvature.  
 By isotopy lemma (see below), the absence of critical point assumed that the manifold is diffeomorphic to the euclidean space $\mathbb{R}^{n}.$

X. Menguy in \cite{11}  and J. Sha and D.Yang in \cite{13} constructed manifolds with nonnegative Ricci curvature and infinite topological type.
Hence  a natural question is under what additional conditions are manifolds with nonnegative Ricci curvature of finite topological type? Are  those manifolds diffeomorphic to the unit sphere or the euclidean space?
Under volume growth, diameter or density of rays conditions, some  results were obtained on the geometry and topology of  open manifolds with nonnegative Ricci curvature. See \cite{2}, \cite{3},\cite{6},\cite{10},\cite{12},\cite{13}, \cite{14},\cite{15}, \cite{17},\cite{18}\ldots

Let $K$ denotes the sectional curvature of $M$ and fix a point $p\in M$. For $r > 0$ let $$k_{p}(r)= \inf_{M\setminus B(p,r)}K$$ where $B(p,r)$ is the open geodesic ball around  with radius $r$ and the infimum is taken over all the sections at points on $M\setminus B(p,r).$
If $(M,g)$  is a complete noncompact Riemannian manifold, we say $M$ has sectional curvature decay at most quadratic if $k_{p}(r)\ge -\frac{C}{r^{\alpha}}$ for some $C > 0,$ $\alpha\in [0,2]$ and all $r> 0.$

In this paper we see the case of manifolds with asymptotically nonnegative Ricci curvature and with sectional curvature decay almost quadratically.

A complete noncompact Riemannian manifold is said to have an asymptotically nonnegative sectional curvature  (Ricci curvature)
if there exists a point $p$, called base point, and a monotne decreasing positive function $\lambda$ such that $\int_{0}^{+\infty}s\lambda(s)ds = b_{0} < +\infty$ and 
 for any point $x$ in $M$ we have \begin{equation*} K(x) \geq -\lambda(d_{p}(x))\; (resp. \;Ric(x)\geq -(n-1)\lambda(d_{p}(x))) 
 \end{equation*} where $d_p$ is the distance to $p$.
Let $B(x,r)$ denote the metric ball of radius $r$ and centre $x$ in $M$ and $B(\overline{x},r)$ denote the similar metric ball in the simply connected noncompact complete manifold with sectional curvature $-\lambda (d_{\overline{p}}(\overline{x}))$ at the point $\overline{x}$ where $d_{p}(x)= d(p,x)$ is the distance from $p$ to $x.$

The  volume comparison theorem proved  in \cite{9}] says that the function $r\mapsto \frac{volB(x,r)}{volB(\overline{x},r)}$ is monotone decreasing. Set $$\alpha_{x} = \lim_{r\rightarrow +\infty}\frac{volB(x,r)}{volB(\overline{x},r)}\; \textrm{ and } \alpha_{M} = \inf_{x\in M}\alpha_{x}.$$ We say $M$ is large volume growth if $\alpha_M > 0.$

In \cite{1}   U. Abresch proved that asymptotically nonnegative sectional curvature have finite topolological type.

Let $R_{p}$ denotes the set of all ray issuing from $p$ and $S(p,r)$ the geodesic ball of radius $r$ and the center $p$. Set $H(p,r) = \max_{x\in S(p,r)}d(x,R_{p}).$  
 By definition,we have $H(p,r)\le r.$
Some results have been obtained by geometers on manifolds with nonnegative Ricci curvature by using the density of the rays. For manifods with quadratic sectional curvature decay, Q. Wang and C. Xia proved that there exists a constant $\delta$ such that if $H(p,r) < \delta r$ then they are diffeomorphic to $\mathbb{R}^{n}.$

In this paper we prove the following theorem:

 \begin{theorem}\label{*}
 
Given $c > 0$ and $\alpha\in[0,2];$ suppose that $M$ is an n-dimensional 
complete noncompact Riemannian manifold with 
$Ricci_{M}(x) \geq -(n-1)\lambda (d_{p}(x))$ and $K(x)\geq -\frac{C}{d_{p}(x)^{\alpha}},\; Crit_{p}\ge r_{0}$ then  there exists 
 a positive constant $\delta _{0} > 0$ such that if $H(p,r)< \delta_{0} r^{\beta/2}$ then $M$ 
 is diffeomorphic to $\mathbb{R}^{n}$ where $\beta = \frac{2}{n} + \alpha (1-\frac{1}{n}).$
 \end{theorem}
 \begin{remark}
 
(i) Theorem\ref{*} is an improvement of theorem1.1 \cite{16} where nonnegative Ricci curvature was assumed and 
 sectional curvature $K_{p}(r)\ge-\frac{C}{(1+ r)^{\alpha}}$.

(ii) For $\alpha = 0$ theorem\ref{*} is a generalisation of lemma 3.1 \cite{18}.
 \end{remark}
 In \cite{16} Q. Wang and C. Xia proved the following theorem (Theorem 1.3)

\begin{theorem}
   Given $\alpha \in [0, 2],$ positive numbers $r_0$ and $C$, and an integer $ n ≥ 2,$ there is an
 $\epsilon = (n, r_0 , C, \beta) > 0$ such that any complete Riemannian $n$-manifold $M$ with Ricci curvature $Ric_M \ge 0,$
$\alpha_M > 0,$ $crit_p \ge r_0$ and $$K(x)\geq -\frac{C}{(1+ d_{p}(x))^{\alpha}},\, \frac{volB(p, r)}{\omega_{n}r^{n}}\le \left(1 +\frac{\epsilon}{r^{(n-2\frac{1}{n})(1-\frac{\alpha}{2})}}\right)\alpha_{M}$$
for some $p \in M$ and all $r \ge r_0$ is diffeomorphic to $\mathbb{R}^{n}.$
\end{theorem}

In this paper we  prove a more general result:

\begin{theorem}\label{a}

Given $c > 0$ and $\alpha\in[0,2];$ suppose that $M$ is an n-dimensional complete noncompact Riemannian manifold with 
$Ricci_{M}(x) \geq -(n-1)\lambda (d_{p}(x))$ and $K(x)\geq -\frac{C}{d_{p}(x)^{\alpha}},\; Crit_{p}\ge r_{0}$ then  there exists a positive constant $\epsilon = \epsilon (C, \alpha, r_{0})$ such that if \begin{equation}\label{c} \frac{volB(p,r)}{volB(\overline{p},r)}\le\left(1+\frac{\epsilon}{r^{(n-2 +\frac{1}{n})(1-\frac{\alpha}{2})}}\right)\alpha_{p} \end{equation} then $M$ is    diffeomorphic to $\mathbb{R}^{n}.$
\end{theorem}
\section{Prelimanaries}
To prove our results we need some lemmas.

 The following one is proved in \cite{8} 
\begin{lemma} (Isotopy Lemma). 

Let $0\le r_{1}\le r_{2}\le\infty.$  If a connected component $C$ of
$B(p, r_2 )\setminus B(p, r_1 )$ is free of critical points of $p$, then $C$ is homeomorphic to $C_1 \times
[r_1 , r_2 ]$, where $C_1$ is a topological submanifold without boundary.
\end{lemma}
If $r_{1} = 0$  and $r_{2} =\infty$ then the homeomorphism becomes diffeomorphism (see for example \cite{7}.)

Let $p$ and $q$ be two points of a complete Riemannian manifold $M.$ The excess function $e_{pq}$ is defined by:
$e_{pq}(x)= d_{p}(x) +d_{q}(x) - d(p,q).$ In \cite{2} U. Abresch and D. Gromoll gave and explicit upper bound of the excess function in manifolds with curvature bounded below. They proved the following lemma:
\begin{lemma}(Proposition 3.1 \cite{2})
Let $M$ be an $n-$dimensional complete Riemannian manifold $(n\ge 3$  and let $\gamma$  be a minimal geodesic joining the base point $p$ and another point $q\in M$ , $x\in M$ is a third point and the excess function  $e_{pq}(x)= d_{p}(x) +d_{q}(x) - d(p,q).$ Suppose  $d(p,q)\ge 2d_p(x)$  and, moreover, that there exists a nonincreasing function $\lambda : [0, +\infty[\rightarrow [0, +\infty[$ such that $b_{0} = \int_{0}^{\infty}r\lambda (r)dr$ converges and $Ric\ge-(n-1)\lambda(d_{p}(x))$ at all points $x\in M.$ Then the height of the triangles can be bounded from below in terms of $d_{p}(x)$ and excess $e_{pq}(x)$. More precisely, \begin{equation}\label{**}s \ge min\left\{\frac{1}{6}d_{p}(x),\frac{d_{p}(x)}{(1+8b_{0})^{1/2}}, C_{0}d_{p}(x)^{1/n}(2e_{pq}(x))^{1-\frac{1}{n}})\right\}\end{equation} where $C_{0} = \frac{4}{17}\frac{n-2}{n-1}(\frac{5}{1+8b_{0}})^{1/n}.$
\end{lemma}

\begin{lemma}[lemma \cite{9}]
\label{d}

Let $(M,g)$ be a complete noncompact Riemannian manifold with asymptotically nonnegative Ricci curvature with  base point $p$
Then for all $x\in M$ and all numbers   $R',R$ with $0 < R' < R$ we have \begin{equation} \frac{volB(x,R)}{volB(x,R')}\le \frac{volB(\overline{x},R)}{vol B(\overline{x},R')}
 \le \left\{ \begin{array}{lcl}
e^{(n-1)b_{0}}\left(\frac{R}{R'}\right)^{n} &if & 0 < R < r=d(p,x)\\
e^{(n-1)b_{0}}\left(\frac{R+r}{R'}\right)^{n}& if& R\geq r
\end{array}\right.\end{equation}
 where $B(\overline{x},s)$  is the ball in
$\overline{M}$ with center $\overline{x}$ and radius $s.$\end{lemma}

\vspace{0.5cm}
Let $\Sigma_{p}$ be a closed subset of $U_{p} = \left\{ u \in T_{p}M, \Vert u\Vert = 1\right\}.$ 

Set $\Sigma_{p}(r)= \left\{ v\in\Sigma_{p}/ \gamma (t) = exp_{p}tv, \gamma \textrm{ is minimal on }
 [0,r]\right\}$
and $$B_{\Sigma_{p}(r)}(p,r) = \left\{ x\in B(p,r) / \exists \gamma : [0,s]\rightarrow M, \gamma (0) = p , \gamma (s) = x \textrm{ and }  \gamma'(0)\in \Sigma_{p}\right\}.$$ Set  $ \Sigma_{p}(\infty)= \cap_{r> 0}\Sigma_{p}(r).$ 

The following two lemmas generalised the above one. 

\begin{lemma}[Lemma3.9 \cite{10}]

Let $(M,g)$ be a Riemannian complete noncompact manifold such that $Ric_{M}\ge -(n-1)\lambda(d_{p}(x))$ 
and $\Sigma_{p}$ be a closed subset of  $U_{p}.$ Then the function 
$r\mapsto\frac{volB_{\Sigma_{p}}(p,r)}{volB(\overline{p},r)}$ is non increasing.\end{lemma}
\begin{lemma}[Lemma 3.10 \cite{10}]\label{b}
Let $(M, g)$ be a Riemannian complete noncompact manifold such that $Ric_{M}\ge -(n-1)\lambda(d_{p}(x))$ and $\Sigma_{p}$ be a closed subset of  $U_{p}.$ Then $\frac{volB_{\Sigma_{p}(r)}(p,r)}{volB(\overline{p},r)}\ge \alpha_{p}.$
\end{lemma}
\section{Proofs}

\textit{Proof of theorem\ref{*}}

To prove the theorem\ref{*}, it suffices to show that $d_p$ has no critical point other than $p$. Let $x$ be a point of $M.$ Set $r = d(p,x);$ $s = d(x, R_{p}).$
Since $R_p$ is closed there exists a ray $\gamma$ issuing from $p$ such that $s = d(x,\gamma)$. 
Set $q = \gamma (t_{0})$ for $t_{0}\ge 2r.$
Let $\sigma_{1}$ and $\sigma_{2}$ be geodesics joining $x$ to $p$ and $q$ respectively.

Set $\tilde{p} = \sigma_{1} (\delta r^{\alpha /2})$ ; $\tilde{q} = \sigma_{2} (\delta r^{\alpha /2})$ with \begin{equation}\label{***}\delta < min\left\{C_{0}^{n},\frac{r_{0}^{1-\beta/2}}{20},\frac{r_{0}^{1-\beta/2}}{\sqrt{1 + 8b_{0}}},
 C_{0}^{n}r_{0}^{1-\beta/2}\right\}.\end{equation}

Consider
the triangle $( x,\tilde{p},\tilde{q});$ if $y$ is a point on this triangle, then
$$d(p,y)\ge d(p,x)-d(x,y)\ge d(p,x)-d(\tilde{p},x)-d(\tilde{p},y)\ge d(p,x)-2\delta r^{\alpha/2}.$$ Since $\beta \ge \alpha.$
we have \begin{equation}d(p,y)\ge d(p,x) - 2\delta r^{\beta/2}\ge r(1-2\delta r^{\frac{\beta}{2}-1})\ge r(1-2\delta r_{0}^{\frac{\beta}{2}-1})\ge r/4.\end{equation}
Hence $y\in M\setminus B(p,r/4)$ and $K_{M}(y)\ge -\frac{4^{\alpha}C}{r^{\alpha}}.$

Thus the triangle $( x,\tilde{p},\tilde{q})\subset M\setminus B(p,\frac{r}{4}).$
Set $\theta =\measuredangle\sigma_{1}'(0),\sigma_{2}'(0).$

Applying the Toponogov's theorem to the triangle $ (x,\tilde{p},\tilde{q})$ we have: 
\begin{equation}\label{p}
 \cosh \left(\frac{2^{\alpha}C^{1/2}}{r^{\alpha/2}}d(\tilde{p},\tilde{q})\right)\le \cosh^{2}\left(\frac{2^{\alpha}C^{1/2}}{r^{\alpha/2}}d(\tilde{p},x)\right)-\sinh^{2}\left(\frac{2^{\alpha}C^{1/2}}{r^{\alpha/2}}d(\tilde{p},x)\right)\cos\theta
\end{equation}
Since $s < \delta r^{\beta/2}$, we deduce from inequaties $(\ref{**})$ and $(\ref{***})$
$$ C_{0}r^{1/n}(2e_{pq}(x))^{1-\frac{1}{n}} < \delta r^{\beta/2},$$ 
hence \begin{equation}
e_{pq}(x)\le \frac{\delta^{n/n-1}}{2C_{0}^{n/n-1}}r^{\alpha/2}  \le \frac{\delta}{2}.r^{\alpha/2}.
\end{equation}
 By triangle inequality, we have
\begin{equation}\left.\begin{array}{rl}d(\tilde{p},\tilde{q})&\ge d(p,q)-d(p,\tilde{p})-d(q,\tilde{q})\\  &\ge d(p,q)-d(p,x)+d(\tilde{p},x)-d(x,q) +d(\tilde{q},x)\\&\ge 2\delta r^{\alpha/2} - e_{pq}(x).\end{array}\right..\end{equation}
Hence
\begin{equation}\label{****}
d(\tilde{p},\tilde{q}) \ge 2\delta r^{\alpha/2}-\frac{\delta}{2}r^{\alpha/2}\ge\frac{3}{2}\delta r^{\alpha/2}
.\end{equation}

From inequalities $(\ref{p})$ and $(\ref{****})$  we deduce 

\begin{displaymath}
 cosh \left(\frac{3}{2}C^{1/2}2^{\alpha}\delta\right)\le cosh^{2} \left(C^{1/2}2^{\alpha}\delta\right)-sinh^{2}\left(C^{1/2}2^{\alpha}\delta\right)cos\theta
.\end{displaymath}
Therefore
$$ sinh^{2}\left(C^{1/2}2^{\alpha}\delta\right)cos\theta\le cosh^{2}\left(C^{1/2}2^{\alpha}\delta\right)-cosh \left(\frac{3}{2}C^{1/2}2^{\alpha}\delta\right)$$ Let $X_{0}$ be the solution of the equation $cosh^{2}2X -cosh3X = 0.$ If $\delta_{0} < \frac{X_{0}}{2^{\alpha-1}}$ then $\theta > \frac{\pi}{2}$ which means that $x$ is not a critical point of $d_{p}$ and the conclusion follows.

\vspace{0.5cm}
\textit{Proof of theorem\ref{a}}

\vspace{0.5cm}
If $y(t)$ denotes the function given by the Jacobi equation $$y''(t)= \lambda (t)y(t)$$ in the simply  connected manifold with sectional curvature $-\lambda (d(\overline{p},\overline{x}))$ at the point $\overline{x}$ then (see \cite{9} ) \begin{equation}\label{q}t\le y(t)\le e^{b_{0}}t\end{equation} and it follows that \begin{equation} \label{e}\omega_{n}r^{n}\le volB(\overline{p},r)\le \omega_{n}e^{(n-1)b_{0}}r^{n}.\end{equation}
In one hand we have:

Let $x\in M$,   $x\ne p$; set $s= d(x,R_{p})$ and $\Sigma_{p}^{c}(\infty)= U_{p}\setminus\Sigma_{p}(\infty).$  Thus
$$ B( x,\frac{s}{2})\subset B_{\Sigma_{p}^{c}(\infty)}(p,r+\frac{s}{2})\setminus B(p,r-\frac{s}{2}).$$
Hence \begin{equation}
 volB( x,\frac{s}{2})\le volB_{\Sigma_{p}^{c}(\infty)}(p,r+\frac{s}{2})-volB(p,r-\frac{s}{2})
 \end{equation}\begin{equation}\le volB_{\Sigma_{p}^{c}(\infty)}(p,r+\frac{s}{2})-volB_{\Sigma_{p}^{c}(\infty)}(p,r-\frac{s}{2})\end{equation}
\begin{equation}\le 
volB_{\Sigma_{p}^{c}(\infty)}(p,r-\frac{r}{2})\left(\frac{volB_{\Sigma_{p}^{c}(\infty)}(p,r+\frac{s}{2})}
{volB_{\Sigma_{p}^{c}(\infty)}(p,r-\frac{r}{2})}-1\right).\end{equation}
We deduce from lemma\ref{b}
 \begin{equation}volB( x,\frac{s}{2})\le volB_{\Sigma_{p}^{c}(\infty)}(p,r-\frac{s}{2})\left(\frac{volB(\overline{p},r+\frac{s}{2})}{volB(\overline{p},r-\frac{s}{2})}-1\right) \end{equation}
\begin{equation}\le \frac{volB_{\Sigma_{p}^{c}(\infty)}(p,r-\frac{s}{2})}{\omega_{n}(r-\frac{s}{2})^{n}}\left(\int_{U_{\overline{p}}}\int_{r-s/2}^{r+s/2}J^{n-1}(t)dt\right)
\end{equation} 
where $J(t)$ denotes the exponential Jacobi in polar coordinates.
Since the function $J/y$ is nonincreasing (see \cite{9}) and using the inequality $(\ref{q}$ we have:

\begin{equation}volB( x,\frac{s}{2})\le \frac{volB_{\Sigma_{p}^{c}(\infty)}(p,r-\frac{s}{2})}{\omega_{n}(r-\frac{s}{2})^{n}} e^{(n-1)b_{0}}\left(\int_{U_{\overline{p}}}\int_{r-s/2}^{r+s/2}t^{n-1}dt\right)\end{equation}
\begin{equation}\le \frac{volB_{\Sigma_{p}^{c}(\infty)}(p,r-\frac{s}{2})}{(r-\frac{s}{2})^{n}} e^{(n-1)b_{0}}\left((r+s/2)^{n}-(r-s/2)^{n}\right)\end{equation} \begin{equation}\le  volB_{\Sigma_{p}^{c}(\infty)}(p,r-\frac{s}{2})e^{(n-1)b_{0}}\left(\left(\frac{r+s/2}{r-s/2}\right)^{n}-1\right)\end{equation}
\begin{equation}\label{r}\le volB_{\Sigma_{p}^{c}(\infty)}(p,r-\frac{s}{2})e^{(n-1)b_{0}}\left((1+\frac{2s}{r})^{n}-1\right)\le volB_{\Sigma_{p}^{c}(\infty)}(p,r-\frac{s}{2})e^{(n-1)b_{0}}.\frac{s}{r}(3^{n}-1).\end{equation}
In other hand we have
$$  volB_{\Sigma_{p}^{c}(\infty)}(p,r-\frac{s}{2})= volB(p,r-s/2)-volB_{\Sigma_{p}(\infty)}(p,r-s/2)$$
By $(\ref{b})$ we have \begin{equation}\label{s}volB_{\Sigma_{p}(\infty)}(p,r-s/2)\ge \alpha_{p}volB(\overline{p},r-s/2).\end{equation}

From  $(\ref{r})$ and $(\ref{s}$ we deduce 
$$volB(x,s/2)\le\left[volB(p,r-s/2)-\alpha_{p}volB(\overline{p},r-s/2)\right].e^{(n-1)b_{0}}.\frac{s}{r}(3^{n}-1).$$ 
By $\ref{c}$ we have 
 \begin{equation}\label{g}volB(x,s/2)\le \frac{\epsilon\alpha_{p}}{(r-s/2)^{(n-2 +\frac{1}{n})(1-\frac{\alpha}{2})}}e^{(n-1)b_{0}}\frac{s}{r}3^{n}volB(\overline{p},r-\frac{s}{2}).
\end{equation}
From $(\ref{e})$ and $(\ref{g})$ we have
\begin{equation}\label{h}volB(x,s/2)\le \epsilon\alpha_{p}e^{2(n-1)b_{0}}s 3^{n}\omega_{n}r^{(n-1)(\frac{1}{n} +\frac{\alpha }{2}(1-\frac{1}{n}))}.\end{equation}
 
We claim that
\begin{equation}\label{i} volB(x,s/2)\ge \frac{\omega_{n}\alpha_{p}}{6^{n}e^{(n-1)b_{0}}}s^{n}.\end{equation}
Indeed we have $B(p,r)\subset B(x,2r),$ and by $(\ref{d})$ we deduce  \begin{equation}\frac{volB(p,r)}{volB(x,s/2)}\le \frac{volB(x,2r)}{volB(x,s/2)}\leq \frac{B(\overline{x},2r)}{volB(\overline{x},s/2)}\end{equation}
\begin{equation}\le e^{(n-1)b_{0}}\left(\frac{2r+s}{s/2}\right)^{n}\le e^{(n-1)b_{0}}6^{n}\left(\frac{r}{s}\right)^{n}.\end{equation}
Thus \begin{equation}\label{f} volB(x,s/2)\ge \frac{s^{n}volB(p,r)}{6^{n}e^{(n-1)b_{0}}r^{n}}.\end{equation}
Hence from $(\ref{e})$, lemma\ref{b} and $(\ref{f})$ the conclusion follows.
 
Thus from $(\ref{h})$ and the  inequality $(\ref{i})$  we have $$s^{n-1}\leq\epsilon 18^{n}e^{3(n-1)b_{0}} r^{(n-1)(\frac{1}{n} +\frac{\alpha }{2}(1-\frac{1}{n}))}$$ which means that $$ s\le \epsilon^{1/(n-1)} 18^{n/(n-1)}e^{3b_{0}} r^{\frac{1}{n} +\frac{\alpha }{2}(1-\frac{1}{n})}.$$
Then it suffices to take $\epsilon < \frac{\delta^{n-1}}{18^{n}e^{3(n-1)b_{0}}}.$

\end{document}